\documentclass[11pt,a4paper]{article}
\usepackage{amsmath,amssymb}
\usepackage{latexsym}
\usepackage{mathrsfs}
\usepackage{amsfonts}
\usepackage{hyperref}
\usepackage{amsthm}
\usepackage{appendix}
\usepackage[square,comma,sort&compress,numbers]{natbib}
\usepackage{color}
\usepackage{subfigure}
\usepackage{graphicx}
\usepackage{float}

\hfuzz=\maxdimen
\theoremstyle{definition}
\newtheorem{lemma}{Lemma}[section]
\newtheorem{definition}[lemma]{Definition}
\newtheorem{theorem}[lemma]{Theorem}
\newtheorem{corollary}[lemma]{Corollary}

\newtheorem{remark}{Remark}

\newtheorem{conjecture}{Conjecture}

\DeclareFixedFont{\Acknowledgment}{OT1}{cmr}{bx}{n}{14pt}
\begin{document}

\title{\bf Combinatorial curvature flows for generalized circle packings on surfaces with boundary}
\author{Xu Xu, Chao Zheng}
\maketitle

\begin{abstract}
In this paper, we investigate the deformation of generalized circle packings on ideally triangulated surfaces with boundary, which is the $(-1,-1,-1)$ type generalized circle packing metric introduced by Guo-Luo \cite{GL2}.
To find hyperbolic metrics on surfaces with totally geodesic boundaries of prescribed lengths, we introduce combinatorial Ricci flow and combinatorial Calabi flow for generalized circle packings on ideally triangulated surfaces with boundary.
Then we prove the longtime existence and global convergence for the solutions of these combinatorial curvature flows,
which provide effective algorithms for finding hyperbolic metrics on surfaces with totally geodesic boundaries of prescribed lengths.
\end{abstract}


\section{Introduction}
Discrete conformal structure on polyhedral surfaces is a discrete analogue of the smooth conformal structure on Riemannian surfaces, which has been extensively studied for decades.
There are lots of researches on discrete conformal structures on closed surfaces. See \cite{BPS,Chow-Luo,GL2,GT,Gu1,Gu2,Gu-Luo-Wu,Guo,Luo1,Luo3,
Luo-Wu,Sun-Wu-Gu-Luo,Thurston,Wu,Wu-Gu-Sun,
Wu-Zhu,Xu1,Xu18, Xu MRL, Xu CAG, XZ TAMS, XZ CVPDE, Zhang-Guo-Zeng-Luo-Yau-Gu,
Stephenson,Guo1,Glickenstein,Bowers} and others for example.
However, there are few researches on discrete conformal structures on surfaces with boundary.
Motivated by Thurston's circle packings on closed surfaces \cite{Thurston}, Guo-Luo \cite{GL2} first introduced some generalized circle packings on surfaces with boundary.
They \cite{GL2} further introduced a generalized combinatorial curvature flow, which is a negative gradient flow of a strictly concave down function. However, no further property of the combinatorial curvature flow is mentioned in \cite{GL2}.
One of the main aims of this paper is to give the longtime behavior of a modification of Guo-Luo's combinatorial curvature flow.
Following Luo's vertex scaling of piecewise linear metrics on closed surfaces \cite{Luo1}, Guo \cite{Guo} introduced a class of hyperbolic discrete conformal structures, also called vertex scaling, on surfaces with boundary.
Guo \cite{Guo} further introduced a combinatorial Yamabe flow for the vertex scaling on surfaces with boundary.
Motivated by the combinatorial Yamabe flow introduced by Guo \cite{Guo}, Li-Xu-Zhou \cite{Li-Xu-Zhou} recently introduced a modified combinatorial Yamabe flow for Guo's vertex scaling on ideally triangulated surfaces with boundary, which generalizes and completes Guo's results \cite{Guo}.
Motivated by Ge \cite{Ge1,Ge2} and Ge-Xu \cite{GX3},
Luo-Xu \cite{Luo-Xu} introduced combinatorial Calabi flow for Guo's vertex scaling on surfaces with boundary and proved its global convergence.
Motivated by the fractional combinatorial Calabi flow introduced by Wu-Xu \cite{Wu-Xu} for discrete conformal structures on closed surfaces, Luo-Xu \cite{Luo-Xu} further introduced fractional combinatorial Calabi flow for Guo's vertex scaling on surfaces with boundary, which unifies and generalizes the combinatorial Yamabe flow and the combinatorial Calabi flow for Guo's vertex scaling on surfaces with boundary. The global convergence of the fractional combinatorial Calabi flow for Guo's vertex scaling on surfaces with boundary is also obtained in \cite{Luo-Xu}.
Recently, Xu \cite{Xu22} introduced a new class of hyperbolic discrete conformal structures on ideally triangulated surfaces with boundary and further introduced the corresponding combinatorial Ricci flow, combinatorial Calabi flow and fractional combinatorial Calabi flow on surfaces with boundary.

In this paper, motivated by the combinatorial Ricci (Yamabe) flow \cite{Guo, Li-Xu-Zhou, Xu22} and the combinatorial Calabi flow \cite{Luo-Xu, Xu22}  on surfaces with boundary, we introduce the combinatorial Ricci flow and the combinatorial Calabi flow for generalized circle packings on ideally triangulated surfaces with boundary, which is the $(-1,-1,-1)$ type generalized circle packing introduced by Guo-Luo \cite{GL2}.
We further prove the longtime existence and global convergence for the solutions of the combinatorial Ricci flow and the combinatorial Calabi flow on  ideally triangulated surfaces with boundary.

Suppose $\Sigma$ is a compact surface with boundary $\partial\Sigma$ consisting of $N$ connected components, which are topologically circles.
Let $\widetilde{\Sigma}$ be the compact surface obtained by coning off each boundary component of $\Sigma$ to be a point,
thus there are exactly $N$ cone points $\{v_1,...,v_N\}$ in $\widetilde{\Sigma}$ so that $\widetilde{\Sigma}-\{v_1,...,v_N\}$ is homeomorphic to $\Sigma-\partial\Sigma$.
An ideal triangulation $\mathcal{T}$ of $\Sigma$ is a triangulation $\widetilde{\mathcal{T}}$ of $\widetilde{\Sigma}$ such that the vertices of the triangulation are exactly the cone points $\{v_1,...,v_N\}$.
The ideal edges and ideal faces of $\Sigma$ in the ideal triangulation $\mathcal{T}$ are defined to be the intersection $\widetilde{E}\bigcap \Sigma$ and $\widetilde{F}\bigcap \Sigma$, where $\widetilde{E}$ and $\widetilde{F}$ are the sets of edges and faces in the triangulation $\widetilde{\mathcal{T}}$ of $\widetilde{\Sigma}$.
The intersection of an ideal face and $\partial\Sigma$ are called boundary arcs.
For simplicity, we denote the connected components of the boundary as $B=\{1, 2,...,N\}$, and the sets of ideal edges and ideal faces as $E$ and $F$ respectively.
The ideal edge between two adjacent boundary components $i,j\in B$ is denoted by $\{ij\}$
and the ideal face adjacent to boundary components $i,j,k\in B$ is denoted by $\{ijk\}$.

The edge length function associated to $\mathcal{T}$ is a vector $l: E\rightarrow (0,+\infty)$ assigning each ideal edge $\{ij\}$ a positive number $l_{ij}$.
For an ideal face $\{ijk\}$ adjacent to boundary components $i,j,k\in B$, there exists a unique hyperbolic right-angled geodesic hexagon whose three non-pairwise adjacent edges having lengths $l_{ij}, l_{jk}, l_{ki}$.
Gluing all such geometric hexagons along the edges in pairs by hyperbolic isometries, one can construct a hyperbolic surface with totally geodesic boundary from the ideal triangulation $\mathcal{T}$.
Conversely, any ideally triangulated hyperbolic surface with totally geodesic boundary $(\Sigma, \mathcal{T})$ produces a function $l: E\rightarrow (0,+\infty)$ with $l_{ij}$ given by the length of the shortest geodesic connecting the boundary components $i,j\in B$.
The edge length function $l: E\rightarrow (0,+\infty)$ is called a discrete hyperbolic metric on $(\Sigma,\mathcal{T})$.
The length $K_i$ of the boundary component $i\in B$ is called the generalized combinatorial curvature of the discrete hyperbolic metric $l:E\rightarrow(0,+\infty)$ at $i\in B$, i.e.,
\begin{equation*}
K_i=\sum_{\{ijk\}\in F}\theta^{jk}_i,
\end{equation*}
where the summation is taken over all the hyperbolic right-angled hexagons with $i$ as a boundary component and $\theta^{jk}_i$ is the length of the boundary arc of the hyperbolic right-angled hexagon $\{ijk\}\in F$ at $i\in B$.

Motivated by Thurston's circle packings on closed surfaces \cite{Thurston},
Guo-Luo \cite{GL2} introduced the following generalized hyperbolic circle packings on surfaces with boundary.
\begin{definition}[\cite{GL2}]\label{def g-cp}
Suppose $(\Sigma,\mathcal{T})$ is an ideally triangulated surface with boundary.
Let $\Phi: E\rightarrow (0,+\infty)$ be the weight on $(\Sigma,\mathcal{T})$.
A generalized circle packing metric of type $(-1,-1,-1)$ is given by a radius function $r: B\rightarrow (0,+\infty)$ so that a discrete hyperbolic metric $l: E\rightarrow (0,+\infty)$ is determined by the radius $r$ and weight $\Phi$ via the following formula
\begin{equation}\label{cosine law}
\cosh l_{ts}=\cosh \Phi_{ts}\sinh r_t\sinh r_s-\cosh r_t\cosh r_s,\ \forall\{st\}\in E.
\end{equation}
\end{definition}

\begin{figure}[!ht]
\centering
\includegraphics[scale=1]{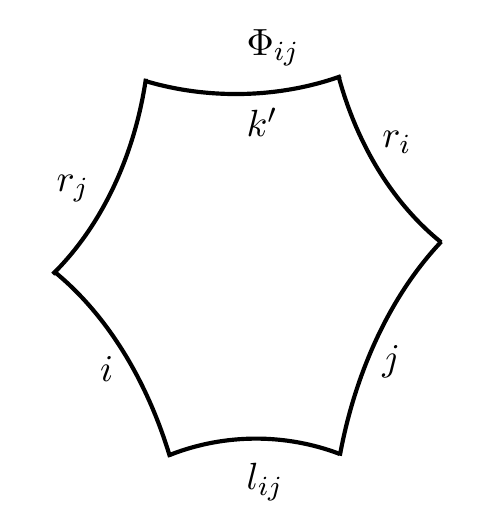}
\caption{Hyperbolic righted-angled hexagon $\{ijk'\}$}
\label{figure1}
\end{figure}

Specially, one can construct a hyperbolic righted-angled hexagon $\{ijk'\}$ such that the edges $\{jk'\}$, $\{ik'\}$ have lengths $r_i,\ r_j$ respectively and the length of the hyperbolic arc at $k'\in B$ is $\Phi_{ij}$.
Let $l_{ij}$ be the length of edge $\{ij\}$ in the hyperbolic righted-angled hexagon $\{ijk'\}$,
then $l_{ij}$ is a function of $r_i,\ r_j$ and fixed $\Phi_{ij}$ defined by (\ref{cosine law}) by the cosine law for hyperbolic right-angled hexagons.
Please refer to Figure \ref{figure1}.
Indeed, to make sure the edge length $l_{ij}$ exists, there are some restrictions on $r_i,\ r_j$.
We need $r_i,\ r_j\in \mathcal{M}(\Phi_{ij})$, which is the admissible space of generalized circle packings $r$ defined by (\ref{admissible space 4}).
One can construct the edge lengths $l_{jk}$, $l_{ki}$ similarly.
A basic fact in hyperbolic geometry is that given any three positive numbers, there exists a unique hyperbolic righted-angled hexagon up to hyperbolic isometry with the lengths of three non-adjacent edges given by the three positive numbers \cite{Ratcliffe}.
Therefore, the edge lengths $l_{ij}$, $l_{jk}$, $l_{ki}$ also determine a hyperbolic righted-angled hexagon $\{ijk\}\in F$.

Motivated by the combinatorial Ricci (Yamabe) flow \cite{Guo, Li-Xu-Zhou, Xu22} and the combinatorial Calabi flow \cite{Luo-Xu, Xu22} on surfaces with boundary, we introduce the following combinatorial Ricci flow and combinatorial Calabi flow for the generalized circle packings in Definition \ref{def g-cp}.

\begin{definition}
Suppose $(\Sigma,\mathcal{T})$ is an ideally triangulated surface with boundary.
Let $\Phi: E\rightarrow (0,+\infty)$ be the weight on $(\Sigma,\mathcal{T})$.
And $\overline{K}\in (0,+\infty)^N$ is a given function defined on $B=\{1,2,...,N\}$.
The combinatorial Ricci flow for the generalized circle packings in Definition \ref{def g-cp} is defined to be
\begin{eqnarray}\label{Ricci flow formula}
\begin{cases}
\frac{du_i}{dt}=K_i-\overline{K}_i,\\
u_i(0)=u_0,
\end{cases}
\end{eqnarray}
where $u_i=\log\tanh\frac{r_i}{2}$ is called a discrete conformal factor.
The combinatorial Calabi flow for the generalized circle packings in Definition \ref{def g-cp} is defined to be
\begin{eqnarray}\label{Calabi flow formula}
\begin{cases}
\frac{du_i}{dt}=-\Delta(K-\overline{K})_i,\\
u_i(0)=u_0,
\end{cases}
\end{eqnarray}
where $\Delta=(\frac{\partial K_i}{\partial u_j})_{N\times N}$ is the discrete Laplace operator.

\end{definition}

\begin{remark}
In the special case of $\overline{K}=0$, the combinatorial Ricci flow (\ref{Ricci flow formula}) differs from Guo-Luo's combinatorial curvature flow by a minus sign \cite{GL2}.
Actually, Guo-Luo's combinatorial curvature flow is a negative gradient flow of a strictly concave down function.
Following this flow, the solution will always diverge.
The combinatorial Ricci flow defined by (\ref{Ricci flow formula}) is a negative gradient flow of a strictly convex function.
\end{remark}

We prove the following result on the combinatorial Ricci flow (\ref{Ricci flow formula}) and the combinatorial Calabi flow (\ref{Calabi flow formula}).

\begin{theorem}\label{main theorem 1}
Suppose $(\Sigma,\mathcal{T})$ is an ideally triangulated surface with boundary.
Let $\Phi: E\rightarrow (0,+\infty)$ be the weight on $(\Sigma,\mathcal{T})$.
For any $\overline{K}\in (0,+\infty)^N$ defined on $B=\{1,2,...,N\}$, the solutions of the combinatorial Ricci flow (\ref{Ricci flow formula}) and the combinatorial Calabi flow (\ref{Calabi flow formula}) exist for all time and converge exponentially fast.
\end{theorem}

By Lemma \ref{matrix negative}, the discrete Laplace operator $\Delta=(\frac{\partial K_i}{\partial u_j})_{N\times N}$ is strictly negative definite on the admissible space $\Omega(\Phi)$ of discrete conformal factors $u$ defined by (\ref{admissible space 1}).
One can define the fractional discrete Laplace operator $\Delta^s$ for any $s\in (-\infty,+\infty)$ as follows.
Recall that if $A$ is a symmetric positive definite $N\times N$ matrix, the matrix $A$ can be written as
\begin{equation*}
A=P^T\cdot \text{diag}\{\lambda_1,...,\lambda_N\}\cdot P,
\end{equation*}
where $P\in O(N)$ is an orthogonal matrix and $\lambda_1,...,\lambda_N$ are positive eigenvalues of the matrix $A$.
For any $s\in (-\infty,+\infty)$, $A^s$ is defined to be
\begin{equation*}
A^s=P^T\cdot \text{diag}\{\lambda^s_1,...,\lambda^s_N\}\cdot P.
\end{equation*}
The $2s$-th order fractional discrete Laplace operator $\Delta^s$ is defined to be
\begin{equation}\label{fractional Laplace}
\Delta^s=-(-\Delta)^s,
\end{equation}
where $\Delta=(\frac{\partial K_i}{\partial u_j})_{N\times N}$ is negative definite by Lemma \ref{matrix negative}.
The fractional discrete Laplace operator $\Delta^s$ is strictly negative definite on the admissible space $\Omega(\Phi)$ by (\ref{fractional Laplace}).
Specially, if $s=0$, the fractional discrete Laplace operator $\Delta^s$ is reduced to the minus identity operator;
if $s=1$, the fractional discrete Laplace operator $\Delta^s$ is reduced to the discrete Laplace operator $\Delta=(\frac{\partial K_i}{\partial u_j})_{N\times N}$.

Motivated by Wu-Xu \cite{Wu-Xu}, we introduce the following fractional combinatorial Calabi flow for generalized circle packings on ideally triangulated surfaces with boundary.

\begin{definition}
Suppose $(\Sigma,\mathcal{T})$ is an ideally triangulated surface with boundary.
Let $\Phi: E\rightarrow (0,+\infty)$ be the weight on $(\Sigma,\mathcal{T})$.
And $\overline{K}\in (0,+\infty)^N$ is a given function defined on $B=\{1,2,...,N\}$.
The fractional combinatorial Calabi flow for the generalized circle packings in Definition \ref{def g-cp} is defined to be
\begin{eqnarray}\label{FCCF formula}
\begin{cases}
\frac{du_i}{dt}=-\Delta^s(K-\overline{K})_i,\\
u_i(0)=u_0,
\end{cases}
\end{eqnarray}
where $\Delta^s$ is the fractional discrete Laplace operator defined by (\ref{fractional Laplace}).
\end{definition}

\begin{remark}
If $s=0$, the fractional combinatorial Calabi flow (\ref{FCCF formula}) is reduced to the combinatorial Ricci flow (\ref{Ricci flow formula}).
If $s=1$, the fractional combinatorial Calabi flow (\ref{FCCF formula}) is reduced to the combinatorial Calabi flow (\ref{Calabi flow formula}).
The fractional combinatorial Calabi flow (\ref{FCCF formula}) further covers the case of $s\neq0,1$.
\end{remark}

\begin{remark}
The notion of the fractional combinatorial Calabi flow was first introduced by Wu-Xu \cite{Wu-Xu} for discrete conformal structures on polyhedral surfaces, which unifies and generalizes Chow-Luo's combinatorial Ricci flow for Thurston's circle packings \cite{Chow-Luo}, Luo's combinatorial Yamabe flow for vertex scaling \cite{Luo1} and the combinatorial Calabi flow for discrete conformal structures on surfaces \cite{Ge1,Ge2,Ge3,GX3,Xu1,Zhu-Xu}.
Motivated by Wu-Xu's fractional combinatorial Calabi flow for discrete conformal structures on polyhedral surfaces,
Luo-Xu \cite{Luo-Xu} introduced fractional combinatorial Calabi flow for Guo's vertex scaling on surfaces with boundary, which unifies and generalizes the combinatorial Yamabe flow \cite{Li-Xu-Zhou} and
the combinatorial Calabi flow \cite{Luo-Xu} on surfaces with boundary.
Xu \cite{Xu22} recently introduced the fractional combinatorial Calabi flow for a new type of
discrete conformal structure on surfaces with boundary.
\end{remark}

We have the following result on the fractional combinatorial Calabi flow (\ref{FCCF formula}).

\begin{theorem}\label{main theorem 2}
Suppose $(\Sigma,\mathcal{T})$ is an ideally triangulated surface with boundary.
Let $\Phi: E\rightarrow (0,+\infty)$ be the weight on $(\Sigma,\mathcal{T})$.
And $\overline{K}\in (0,+\infty)^N$ is a given function defined on $B=\{1,2,...,N\}$.
If the solution $u(t)$ of the fractional combinatorial Calabi flow (\ref{FCCF formula}) converges to $\overline{u}\in \Omega(\Phi)$, then $K(\overline{u})=\overline{K}$.
Furthermore, for any $\overline{K}\in (0,+\infty)^N$ defined on $B=\{1,2,...,N\}$, there exists a constant $\delta >0$ such that if $||K(u(0))-\overline{K}||=
\sqrt{\sum_{i=1}^{N}(K_i(u(0))-\overline{K}_i)^2}<\delta$, then the solution of the fractional combinatorial Calabi flow (\ref{FCCF formula}) exists for all time and converges exponentially fast.
\end{theorem}

The paper is organized as follows.
In Section \ref{section 2}, we give some results derived by Guo-Luo \cite{GL2} on the generalized circle packings
in Definition \ref{def g-cp}.
In Section \ref{section 3}, we study some basic properties of the combinatorial Ricci flow (\ref{Ricci flow formula}) and the combinatorial Calabi flow (\ref{Calabi flow formula}).
In Section \ref{section 4}, we study the longtime behavior of the combinatorial Ricci flow (\ref{Ricci flow formula}) and prove the first part of Theorem \ref{main theorem 1}.
In Section \ref{section 5}, we study the longtime behavior of the combinatorial Calabi flow (\ref{Calabi flow formula}) and prove the second part of Theorem \ref{main theorem 1}.
In Section \ref{section 6}, we give the basic properties of the fractional combinatorial Calabi flow (\ref{FCCF formula}) and prove Theorem \ref{main theorem 2}.
A conjecture about the longtime existence and global convergence for the solution of the fractional combinatorial Calabi flow (\ref{FCCF formula}) on ideally triangulated surfaces with boundary is also proposed in Section \ref{section 6}.\\

\textbf{Acknowledgements}\\[8pt]
The authors would like to thank the referees for their careful reading of the paper and very helpful suggestions.

\section{Preliminaries on generalized circle packings}\label{section 2}
Suppose $(\Sigma,\mathcal{T})$ is an ideally triangulated surface with boundary.
Let $\Phi: E\rightarrow (0,+\infty)$ be the weight on $(\Sigma,\mathcal{T})$.
The admissible space of generalized circle packings $r$ on $(\Sigma,\mathcal{T})$ is defined to be
\begin{equation}\label{admissible space 3}
\mathcal{M}(\Phi)=\bigcap_{\{ij\}\in E}\mathcal{M}(\Phi_{ij}),
\end{equation}
where
\begin{equation}\label{admissible space 4}
\mathcal{M}(\Phi_{ij})
=\{(r_i,r_j)\in (0,+\infty)^2|\ \cosh \Phi_{ij}\sinh r_i\sinh r_j-\cosh r_i\cosh r_j>1\}.
\end{equation}

By a change of variables, i.e., $u_i=\log\tanh\frac{r_i}{2}$, the admissible space $\mathcal{M}(\Phi)$ of generalized circle packings $r$ is transferred to the following admissible space $\Omega(\Phi)$ of discrete conformal factors $u$ on $(\Sigma,\mathcal{T})$
\begin{equation}\label{admissible space 1}
\Omega(\Phi)=\bigcap_{\{ijk\}\in F}\Omega(\Phi_{ij},\Phi_{jk},\Phi_{ki}),
\end{equation}
where
\begin{equation}\label{admissible space 2}
\Omega(\Phi_{ij},\Phi_{jk},\Phi_{ki})=\{(u_i,u_j,u_k)\in (-\infty,0)^3|u_r+u_s>-\Phi_{rs}, \{r,s,t\}=\{i,j,k\} \}.
\end{equation}
The relationship of the admissible space $\mathcal{M}(\Phi)$ and $\Omega(\Phi)$ can be obtained directly from Lemma 4.2 in Guo-Luo \cite{GL2} or derived by calculations in Lemma \ref{Ricci boundary 1}.
It is obvious that the admissible space $\Omega(\Phi_{ij},\Phi_{jk},\Phi_{ki})$ is a convex polyhedron.
As a result, the admissible space $\Omega(\Phi)$ is a convex polyhedron.

Suppose $\{ijk\}\in F$ is a right-angled hyperbolic hexagon adjacent to the boundary components $i,j,k\in B$.
Denote the edge lengths of three nonadjacent edges $\{i j\},\ \{ik\},\ \{jk\}\in E$ as $l_{ij},\ l_{ik},\ l_{jk}$ respectively and the lengths of the hyperbolic arc in the boundary components $i,j,k\in B$ as $\theta_{i}^{jk},\ \theta_{j}^{ik},\ \theta_{k}^{ij}$ respectively.

\begin{lemma}(\cite{GL2})\label{matrix negative}
For any fixed $(\Phi_{ij},\Phi_{jk},\Phi_{ki})\in (0,+\infty)^3$,
the matrix $\frac{\partial (\theta_{i}^{jk}, \theta_{j}^{ik},\theta_{k}^{ij})}{\partial (u_i, u_j, u_k)}$ is symmetric and negative definite on $\Omega(\Phi_{ij},\Phi_{jk},\Phi_{ki})$.
As a result, the matrix $\Delta=\frac{\partial (K_i,..., K_N)}{\partial(u_i,...,u_N)}$ is symmetric and negative definite on $\Omega(\Phi)$.
\end{lemma}

For the generalized circle packing in Definition \ref{def g-cp}, Guo-Luo \cite{GL2} further prove the following rigidity and existence theorem.

\begin{theorem}(\cite{GL2})\label{rigidity theorem}
Suppose $(\Sigma,\mathcal{T})$ is an ideally triangulated surface with boundary.
Let $\Phi: E\rightarrow (0,+\infty)$ be the weight on $(\Sigma,\mathcal{T})$.
Then the generalized circle packing $r\in \mathcal{M}(\Phi)$ is determined by its generalized discrete curvature $K: B\rightarrow (0,+\infty)^N$.
In particular, the map $K$ is a smooth embedding.
Furthermore, the image of $K$ is $(0,+\infty)^N$.
\end{theorem}

\section{Combinatorial curvature flows on surfaces with boundary}\label{section 3}
In this section, we first study some basic properties of the combinatorial Ricci flow (\ref{Ricci flow formula}) and the combinatorial Calabi flow (\ref{Calabi flow formula}).
Then we give some results on the longtime existence and convergence for the solutions of the combinatorial Ricci flow (\ref{Ricci flow formula}) and the combinatorial Calabi flow (\ref{Calabi flow formula}) for initial
data with small energy.

\subsection{Two important energy functions}

Consider the following two energy functions
\begin{equation}\label{energy function 1}
\mathcal{E}(u)
=-\int_0^{u}\sum_{i=1}^{N}(K_i-\overline{K}_i)du_i
\end{equation}
and
\begin{equation}\label{energy function 2}
\mathcal{C}(u)
=\frac{1}{2}||K-\overline{K}||^2
=\frac{1}{2}\sum_{i=1}^{N}(K_i-\overline{K}_i)^2.
\end{equation}

Combining the fact that the admissible space $\Omega(\Phi)$ of discrete conformal factors $u$ defined by (\ref{admissible space 1}) is a convex polyhedron and Lemma \ref{matrix negative},
the function $\mathcal{E}(u)$ is a well-defined strictly convex function defined on the admissible space $\Omega(\Phi)$.

\begin{lemma}\label{CRF lemma}
The combinatorial Ricci flow (\ref{Ricci flow formula}) is a negative gradient flow of the function $\mathcal{E}(u)$ defined by (\ref{energy function 1}).
As a result, the function $\mathcal{E}(u)$ defined by (\ref{energy function 1}) is decreasing along the combinatorial Ricci flow (\ref{Ricci flow formula}).
Furthermore, the combinatorial Calabi energy $\mathcal{C}(u)$ defined by (\ref{energy function 2}) is decreasing along the combinatorial Ricci flow (\ref{Ricci flow formula}).
\end{lemma}
\proof
By direct calculations, we have
\begin{equation*}
\nabla_{u_i}\mathcal{E}(u)=-(K_i-\overline{K}_i),
\end{equation*}
which implies the combinatorial Ricci flow (\ref{Ricci flow formula}) is a negative gradient flow of the function $\mathcal{E}(u)$. Furthermore,
\begin{equation*}
\frac{d\mathcal{E}(u(t))}{dt}
=\sum_{i=1}^{N}\frac{\partial \mathcal{E}}{\partial u_i}\frac{d u_i}{dt}
=-\sum_{i=1}^{N}(K_i-\overline{K}_i)^2\leq0,
\end{equation*}
which implies that the function $\mathcal{E}(u)$ is decreasing along the combinatorial Ricci flow (\ref{Ricci flow formula}). Similarly,
\begin{equation*}
\frac{d\mathcal{C}(u(t))}{dt}
=\sum_{i,j=1}^{N}\frac{\partial \mathcal{C}}{\partial K_i}\frac{\partial K_i}{\partial u_j}\frac{d u_j}{dt}
=\sum_{i,j=1}^{N}(K_i-\overline{K}_i)\frac{\partial K_i}{\partial u_j}(K_j-\overline{K}_j)
=(K-\overline{K})^T\Delta(K-\overline{K})\leq0
\end{equation*}
by Lemma \ref{matrix negative}, the right side of which is strictly negative unless $K=\overline{K}$.
\qed

\begin{lemma}\label{CCF lemma}
The combinatorial Calabi flow (\ref{Calabi flow formula}) is a negative gradient flow of the combinatorial Calabi energy $\mathcal{C}(u)$ defined by (\ref{energy function 2}).
As a result, the combinatorial Calabi energy $\mathcal{C}(u)$ defined by (\ref{energy function 2}) is decreasing along the combinatorial Calabi flow (\ref{Calabi flow formula}).
Furthermore, the energy function $\mathcal{E}(u)$ defined by (\ref{energy function 1}) is decreasing along the combinatorial Calabi flow (\ref{Calabi flow formula}).
\end{lemma}
\proof
By direct calculations, we have
\begin{equation*}
\nabla_{u_i}\mathcal{C}(u)=\sum_{i=1}^N\frac{\partial K_j}{\partial u_i}(K_j-\overline{K}_j)
=\Delta(K-\overline{K})_i=-\frac{du_i}{dt},
\end{equation*}
which implies that the combinatorial Calabi flow (\ref{Calabi flow formula}) is a negative gradient flow of the combinatorial Calabi energy $\mathcal{C}(u)$.
Furthermore,
\begin{equation*}
\frac{d\mathcal{C}(u(t))}{dt}
=\sum_{i=1}^{N}\frac{\partial \mathcal{C}}{\partial u_i}\frac{d u_i}{dt}
=-\sum_{i=1}^{N}(\Delta(K-\overline{K})_i)^2\leq0,
\end{equation*}
the right side of which is strictly negative unless $K=\overline{K}$.
Similarly,
\begin{equation*}
\frac{d\mathcal{E}(u(t))}{dt}
=\sum_{i=1}^{N}\frac{\partial \mathcal{E}}{\partial u_i}\frac{d u_i}{dt}
=\sum_{i=1}^{N}(K-\overline{K})_i\Delta(K-\overline{K})_i
=(K-\overline{K})^T\Delta(K-\overline{K})\leq0,
\end{equation*}
where the last inequality follows from the strictly negative definiteness of the discrete Laplace operator $\Delta=(\frac{\partial K_i}{\partial u_j})_{N\times N}$ in Lemma \ref{matrix negative}.
\qed

\subsection{Local convergence of combinatorial curvature flows}
Note that the combinatorial Ricci flow (\ref{Ricci flow formula}) and the combinatorial Calabi flow (\ref{Calabi flow formula}) are ODE systems with smooth coefficients.
Therefore, the solutions always exist locally around the initial time $t=0$ by the standard ODE theory.
We further have the following result on the longtime existence and convergence for the solutions of the combinatorial Ricci flow (\ref{Ricci flow formula})
and the combinatorial Calabi flow (\ref{Calabi flow formula}) with small initial energy.

\begin{theorem}\label{local converge}
Suppose $(\Sigma,\mathcal{T})$ is an ideally triangulated surface with boundary.
Let $\Phi: E\rightarrow (0,+\infty)$ be the weight on $(\Sigma,\mathcal{T})$. And $\overline{K}\in (0,+\infty)^N$ is a given function defined on $B=\{1,2,...,N\}$.
If the solution $u(t)$ of the combinatorial Ricci flow (\ref{Ricci flow formula}) or the combinatorial Calabi flow (\ref{Calabi flow formula}) converges to $\overline{u}\in \Omega(\Phi)$, then $K(\overline{u})=\overline{K}$.
Furthermore, for any $\overline{K}\in (0,+\infty)^N$ defined on $B=\{1,2,...,N\}$, there exists a constant $\delta >0$ such that if $||K(u(0))-\overline{K}||=
\sqrt{\sum_{i=1}^{N}(K_i(u(0))-\overline{K}_i)^2}<\delta$, then the solution of the combinatorial Ricci flow (\ref{Ricci flow formula}) (the combinatorial Calabi flow (\ref{Calabi flow formula}) respectively) exists for all time and converges exponentially fast.
\end{theorem}
\proof
Suppose $u(t)$ is a solution of the combinatorial Ricci flow (\ref{Ricci flow formula}).
If $\overline{u}:=u(+\infty)=\lim_{t\rightarrow +\infty} u(t)$ exists in $\Omega(\Phi)$, then $K(\overline{u})=\lim_{t\rightarrow +\infty}K(u(t))$ exists by the $C^1$-smoothness of $K$.
Furthermore, there exists a sequence $\xi_n\in(n,n+1)$ such that as $n\rightarrow +\infty$,
$$u_i(n+1)-u_i(n)=u'_i(\xi_n)=K_i(u(\xi_n))-\overline{K}_i\rightarrow 0,$$
which implies $K(\overline{u})=\overline{K}$.
Similarly, if the solution $u(t)$ of the combinatorial Calabi flow (\ref{Calabi flow formula}) converges, then
$K(\overline{u})=\lim_{t\rightarrow +\infty}K(u(t))$ exists by the $C^1$-smoothness of $K$.
Furthermore, there exists a sequence $\xi_n\in(n,n+1)$ such that as $n\rightarrow +\infty$,
\begin{equation*}
u_i(n+1)-u_i(n)=u'_i(\xi_n)
=\Delta(K(u(\xi_n))-\overline{K})_i\rightarrow 0,
\end{equation*}
which implies $K(\overline{u})=\overline{K}$ by the strictly negative definiteness of the discrete Laplace operator $\Delta$ in Lemma \ref{matrix negative}.

For the combinatorial Ricci flow (\ref{Ricci flow formula}), set $\Gamma(u)=K-\overline{K}$.
Then $D\Gamma|_{u=\overline{u}}=\Delta$ is negative definite by Lemma \ref{matrix negative}, which implies $\overline{u}$ is a local attractor of the combinatorial Ricci flow (\ref{Ricci flow formula}).
Then the conclusion follows from Lyapunov Stability Theorem (\cite{Pontryagin}, Chapter 5).
Similarly, for the combinatorial Calabi flow (\ref{Calabi flow formula}), set $\Gamma(u)=-\Delta(K-\overline{K})$.
Then $D\Gamma|_{u=\overline{u}}=-\Delta^2$ is negative definite by Lemma \ref{matrix negative}, which implies $\overline{u}$ is a local attractor of the combinatorial Calabi flow (\ref{Calabi flow formula}).
Then the conclusion follows from Lyapunov Stability Theorem (\cite{Pontryagin}, Chapter 5).
\qed

Theorem \ref{local converge} gives the longtime existence and convergence for the solutions of the combinatorial Ricci flow (\ref{Ricci flow formula}) and the combinatorial Calabi flow (\ref{Calabi flow formula}) for initial value with small energy.
For general initial value, we can further prove that the solutions of the combinatorial Ricci flow (\ref{Ricci flow formula}) and the combinatorial Calabi flow (\ref{Calabi flow formula}) can not reach the boundary of the admissible space, i.e., the combinatorial Ricci flow (\ref{Ricci flow formula}) and the combinatorial Calabi flow (\ref{Calabi flow formula}) can not develop singularities.

\section{The longtime behavior of combinatorial Ricci flow}\label{section 4}
In this section, we prove the longtime existence and convergence of the solution of the combinatorial Ricci flow (\ref{Ricci flow formula}) for general initial value, which is the first part of Theorem \ref{main theorem 1}.

\begin{lemma}\label{u bounded 1}
Suppose $(\Sigma,\mathcal{T})$ is an ideally triangulated surface with boundary.
Let $\Phi: E\rightarrow (0,+\infty)$ be the weight on $(\Sigma,\mathcal{T})$.
For any $\overline{K}\in (0,+\infty)^N$ defined on $B=\{1,2,...,N\}$,
the solution $u(t)$ of the combinatorial Ricci flow (\ref{Ricci flow formula}) stays in a bounded subset of $(-\infty,0)^N$.
\end{lemma}
\proof
By Theorem \ref{rigidity theorem}, there exists $\overline{u}\in \Omega(\Phi)$ such that $K(\overline{u})=\overline{K}$, which implies $\nabla\mathcal{E}(\overline{u})
=-(K-\overline{K})|_{u=\overline{u}}=0.$
Note that $\mathcal{E}(u)$ is a strictly convex function on $\Omega(\Phi)\subseteq(-\infty,0)^N$, we have
\begin{equation}\label{E(u) infty}
\lim_{u\rightarrow -\infty}\mathcal{E}(u)=+\infty.
\end{equation}
Note that $\mathcal{E}(u)$ is decreasing along the combinatorial Ricci flow (\ref{Ricci flow formula}) by Lemma \ref{CRF lemma}, then
\begin{equation}\label{E(u) bounded}
\mathcal{E}(u(t))\leq \mathcal{E}(u(0)).
\end{equation}
Combining (\ref{E(u) infty}) and (\ref{E(u) bounded}), the solution $u(t)$ of the combinatorial Ricci flow (\ref{Ricci flow formula}) stays in a bounded subset of $(-\infty,0)^N$.
\qed

To prove the longtime existence of the solution $u(t)$ of the combinatorial Ricci flow (\ref{Ricci flow formula}),
we just need to prove that the solution $u(t)$ of the combinatorial Ricci flow (\ref{Ricci flow formula}) stays in a compact subset of the admissible space $\Omega(\Phi)$, i.e., the solution $u(t)$ can not reach the boundary of the admissible space $\Omega(\Phi)$.
By (\ref{admissible space 1}) and (\ref{admissible space 2}), the boundary of the admissible space $\Omega(\Phi)$ in $[-\infty,0]^N$ consists of the following three parts
\begin{description}
\item[(1)]
$\partial_{\infty}\Omega(\Phi)=\{u\in[-\infty,0]^N|\ \text{there exists at least}\ i\in B\ \text{such that}\  u_i=-\infty \}$,
\item[(2)]
$\partial_0\Omega(\Phi)=\{u\in[-\infty,0]^N|\ \text{there exists at least}\ i\in B\ \text{such that}\  u_i=0 \}$,
\item[(3)]
$\partial_{l}\Omega(\Phi)
=\bigcup_{\{ij\}\in E}\partial_{ij}\Omega(\Phi)
=\bigcup_{\{ij\}\in E}\{u\in [-\infty,0]^N|u_i+u_j=-\Phi_{ij}\}$.
\end{description}

Lemma \ref{u bounded 1} shows that the solution $u(t)$ of the combinatorial Ricci flow (\ref{Ricci flow formula}) can  not reach the boundary $\partial_{\infty}\Omega(\Phi)$.

\begin{lemma}\label{Ricci boundary 1}
The solution $u(t)$ of the combinatorial Ricci flow (\ref{Ricci flow formula}) can not reach the boundary $\partial_{l}\Omega(\Phi)$.
\end{lemma}
\proof
Suppose $\{ijk\}\in F$ is a right-angled hyperbolic hexagon adjacent to the boundary components $i,j,k\in B$.
For any fixed $(\Phi_{ij},\Phi_{jk},\Phi_{ki})\in (0,+\infty)^3$,
if $u_i+u_j=-\Phi_{ij}$, by $u_i=\log\tanh\frac{r_i}{2}$,
we have
\begin{equation*}
\log\tanh\frac{r_i}{2}+\log\tanh\frac{r_j}{2}=-\Phi_{ij},
\end{equation*}
which is equivalent to
\begin{equation*}
\tanh\frac{r_i}{2}\tanh\frac{r_j}{2}=e^{-\Phi_{ij}}.
\end{equation*}
Hence
\begin{equation}\label{formula 5}
\frac{1}{2}(\tanh\frac{r_i}{2}\tanh\frac{r_j}{2}
+\frac{1}{\tanh\frac{r_i}{2}\tanh\frac{r_j}{2}})
=\frac{1}{2}(e^{-\Phi_{ij}}+e^{\Phi_{ij}})=\cosh \Phi_{ij}.
\end{equation}
Note that
\begin{equation*}
\frac{1}{2}(\tanh\frac{r_i}{2}\tanh\frac{r_j}{2}
+\frac{1}{\tanh\frac{r_i}{2}\tanh\frac{r_j}{2}})
=\frac{1+\cosh r_i\cosh r_j}{\sinh r_i\sinh r_j},
\end{equation*}
and we have
\begin{equation*}
\cosh \Phi_{ij}=\frac{\cosh l_{ij}+\cosh r_i\cosh r_j}{\sinh r_i\sinh r_j}
\end{equation*}
by (\ref{cosine law}).
This implies $l_{ij}=0$ by (\ref{formula 5}).
Therefore, if $u_i+u_j=-\Phi_{ij}$, then the edge length $l_{ij}=0$.
Similarly, we can obtain the edge length $l_{ij}>0$ is equivalent to $u_i+u_j>-\Phi_{ij}$,
which shows the equivalence of the admissible space $\mathcal{M}(\Phi)$ and $\Omega(\Phi)$.

By the cosine law for a hyperbolic right-angled hexagon, we have
\begin{equation*}
\cosh \theta^{jk}_i=\frac{\cosh l_{jk}+\cosh l_{ij}\cosh l_{ik}}{\sinh l_{ij}\sinh l_{ik}}>\frac{\cosh l_{ij}\cosh l_{ik}}{\sinh l_{ij}\sinh l_{ik}}>\frac{\cosh l_{ij}}{\sinh l_{ij}},
\end{equation*}
which implies that $\theta^{jk}_i\rightarrow +\infty$ uniformly as $l_{ij}\rightarrow 0$.
Note that $K_i=\sum_{\{ijk\}\in F}\theta^{jk}_i$ and the combinatorial Calabi energy $\mathcal{C}(u)$ is decreasing along the combinatorial Ricci flow (\ref{Ricci flow formula}) by Lemma \ref{CRF lemma}, then
\begin{equation*}
\mathcal{C}(u(t))\leq \mathcal{C}(u(0)),
\end{equation*}
which implies that for any $i\in B$,
$|K_i-\overline{K}_i|$ is bounded along the combinatorial Ricci flow (\ref{Ricci flow formula}).
Since $\overline{K}$ is a given function, then
$K$ is bounded from above along the combinatorial Ricci flow (\ref{Ricci flow formula}).
As a result, $\theta^{jk}_i$ can not tend to $+\infty$ along the combinatorial Ricci flow (\ref{Ricci flow formula}), which implies that the edge length $l_{ij}$ has a positive lower bound along the combinatorial Ricci flow (\ref{Ricci flow formula}).
This is equivalent to $u_i+u_j>-\Phi_{ij}$.
Therefore, the solution $u(t)$ of the combinatorial Ricci flow (\ref{Ricci flow formula}) can not reach the boundary $\partial_{l}\Omega(\Phi)$.
\qed

To prove the solution $u(t)$ of the combinatorial Ricci flow (\ref{Ricci flow formula}) can not reach the boundary $\partial_{0}\Omega(\Phi)$, we need the following lemma, which is derived by Guo-Luo \cite{GL2}.

\begin{lemma}(\cite{GL2}, Lemma 4.6)\label{Guo-Luo Lemma}
Suppose $\{ijk\}\in F$ is a right-angled hyperbolic hexagon adjacent to the boundary components $i,j,k\in B$.
For any fixed $(\Phi_{ij},\Phi_{jk},\Phi_{ki})\in (0,+\infty)^3$, then
\begin{equation*}
\lim_{r_k\rightarrow+\infty}\theta^{ij}_k(r_i,r_j,r_k)=0
\end{equation*}
and the converge is uniform.
\end{lemma}

\begin{lemma}\label{Ricci boundary 2}
The solution $u(t)$ of the combinatorial Ricci flow (\ref{Ricci flow formula}) can not reach the boundary $\partial_0\Omega(\Phi)$.
\end{lemma}
\proof
The trick is taken from Ge-Xu \cite{Ge-Xu 17}.
Suppose that $\lim_{t\rightarrow T}u_i(t)=0$ for $T\in(0,+\infty]$.
By $u_i=\log\tanh\frac{r_i}{2}$, we have $\lim_{t\rightarrow T}r_i(t)=+\infty$.
By Lemma \ref{Guo-Luo Lemma}, $\theta^{jk}_{i}\rightarrow 0$ uniformly as $r_i\rightarrow+\infty$.
Therefore, there exists $c\in (-\infty,0)$ such that whenever $u_i(t)>c$,  the length of the hyperbolic arc $\theta^{jk}_{i}$ is smaller than $\epsilon>0$,
which implies $K_i<\overline{K}_i$.
Choose a time $t_0\in (0,T)$ such that $u_i(t_0)>c$, this can be done because $\lim_{t\rightarrow T}u_i(t)=0$.
Set $a=\inf\{t<t_0|u_i(s)>c, \forall s\in (t,t_0]\}$, then $u_i(a)=c$.
Note that $\frac{du_i}{dt}=K_i-\overline{K}_i<0$ on $(a,t_0]$, we have $u_i(t_0)<u_i(a)=c$, which contradicts $u_i(t_0)>c$.
Therefore, $u_i(t)$ is uniformly bounded from above in $(-\infty,0)$, i.e., the solution $u(t)$ of the combinatorial Ricci flow (\ref{Ricci flow formula}) can not reach the boundary $\partial_{0}\Omega(\Phi)$.
\qed

As a direct corollary of Lemma \ref{u bounded 1}, Lemma \ref{Ricci boundary 1} and Lemma \ref{Ricci boundary 2}, we have the following result on the solution of the combinatorial Ricci flow (\ref{Ricci flow formula}), which implies the longtime existence of the solution of the combinatorial Ricci flow (\ref{Ricci flow formula}).

\begin{corollary}\label{compact subset 1}
Suppose $(\Sigma,\mathcal{T})$ is an ideally triangulated surface with boundary.
Let $\Phi: E\rightarrow (0,+\infty)$ be the weight on $(\Sigma,\mathcal{T})$.
For any $\overline{K}\in (0,+\infty)^N$ defined on $B=\{1,2,...,N\}$,
the solution $u(t)$ of the combinatorial Ricci flow (\ref{Ricci flow formula}) stays in a compact subset of the admissible space $\Omega(\Phi)$.
As a result, the solution $u(t)$ of the combinatorial Ricci flow (\ref{Ricci flow formula}) exists for all time.
\end{corollary}

The following theorem shows the global convergence of the solution of the combinatorial Ricci flow (\ref{Ricci flow formula}), which is the first part of Theorem \ref{main theorem 1}.

\begin{theorem}
Suppose $(\Sigma,\mathcal{T})$ is an ideally triangulated surface with boundary.
Let $\Phi: E\rightarrow (0,+\infty)$ be the weight on $(\Sigma,\mathcal{T})$.
For any $\overline{K}\in (0,+\infty)^N$ defined on $B=\{1,2,...,N\}$, the solution of the combinatorial Ricci flow (\ref{Ricci flow formula}) converges exponentially fast.
\end{theorem}
\proof
By Theorem \ref{rigidity theorem}, there exists a discrete conformal factor $\overline{u}$ such that $\overline{K}=K(\overline{u})$.
By Lemma \ref{CRF lemma} and Corollary \ref{compact subset 1}, we have $\lim_{t\rightarrow +\infty}\mathcal{E}(u(t))$ exists.
Therefore, there exists a sequence $\xi_n\in(n,n+1)$ such that as $n\rightarrow +\infty$,
\begin{equation*}
\mathcal{E}(u(n+1))-\mathcal{E}(u(n))
=(\mathcal{E}(u(t))'|_{\xi_n}=\nabla \mathcal{E}\cdot\frac{du_i}{dt}|_{\xi_n}
=-\sum_{i=1}^{N}(K_i(u(\xi_n))-\overline{K}_i)^2
\rightarrow 0.
\end{equation*}
Then $\lim_{n\rightarrow +\infty}K_i(u(\xi_n))=\overline{K}_i=K_i(\overline{u})$ for all $i\in B$.
Since $\{u(t)\}\subset\subset \Omega(\Phi)$ by Corollary \ref{compact subset 1}, there exists $u^*\in \Omega(\Phi)$ and a subsequence of $\{u(\xi_n)\}$, still denoted as $\{u(\xi_n)\}$ for simplicity,
such that $\lim_{n\rightarrow \infty}u(\xi_n)=u^*$, which implies
$K_i(u^*)=\lim_{n\rightarrow +\infty}K_i(u(\xi_n))=K_i(\overline{u})$.
This further implies $\overline{u}=u^*$ by Theorem \ref{rigidity theorem}.
Therefore, $\lim_{n\rightarrow \infty}u(\xi_n)=\overline{u}$.

Set $\Gamma(u)=K-\overline{K}$, then the matrix $D\Gamma|_{u=\overline{u}}=\Delta<0$ by Lemma \ref{matrix negative}, i.e. $D\Gamma|_{u=\overline{u}}$ has $N$ negative eigenvalues, which implies that $\overline{u}$ is a local attractor of the combinatorial Ricci flow (\ref{Ricci flow formula}). Then the conclusion follows from Lyapunov Stability Theorem (\cite{Pontryagin}, Chapter 5).
\qed

\section{The longtime behavior of combinatorial Calabi flow}\label{section 5}
In this section, we prove the longtime existence and convergence of the solution of the combinatorial Calabi flow (\ref{Calabi flow formula}) for general initial value, which is the second part of Theorem \ref{main theorem 1}.

The following two lemmas are paralleling to Lemma \ref{u bounded 1} and Lemma \ref{Ricci boundary 1}.
As the proofs are all the same, we omit them here.
\begin{lemma}\label{u(t) bounded 2}
Suppose $(\Sigma,\mathcal{T})$ is an ideally triangulated surface with boundary.
Let $\Phi: E\rightarrow (0,+\infty)$ be the weight on $(\Sigma,\mathcal{T})$.
For any $\overline{K}\in (0,+\infty)^N$ defined on $B=\{1,2,...,N\}$,
the solution $u(t)$ of the combinatorial Calabi flow (\ref{Calabi flow formula}) stays in a bounded subset of $(-\infty,0)^N$,
which implies the solution $u(t)$ of the combinatorial Calabi flow (\ref{Calabi flow formula}) can not reach the boundary $\partial_{\infty}\Omega(\Phi)$.
\end{lemma}

\begin{lemma}\label{CCF boundary 1}
The solution $u(t)$ of the combinatorial Calabi flow (\ref{Calabi flow formula}) can not reach the boundary $\partial_{l}\Omega(\Phi)$.
\end{lemma}

To prove the solution $u(t)$ of the combinatorial Calabi flow (\ref{Calabi flow formula}) can not reach the boundary $\partial_{0}\Omega(\Phi)$, we need the following lemma.

\begin{lemma}\label{main lemma}
Suppose $\{ijk\}\in F$ is a right-angled hyperbolic hexagon adjacent to the boundary components $i,j,k\in B$.
Given any $(\Phi_{ij},\Phi_{jk},\Phi_{ki})\in (0,+\infty)^3$, set $l_p=l_{st},\ \theta_p=\theta^{st}_p,\ \{p,s,t\}=\{i,j,k\}$ for simplicity.
Then for any $C\in (-\infty,+\infty)$, there exists a constant $M=M(C)>0$ such that if $r_i\geq M$, then
\begin{equation}\label{main formula}
\bigg|\frac{\partial \theta_i}{\partial u_i}\bigg|
>C\bigg(\bigg|\frac{\partial \theta_i}{\partial u_j}\bigg|+\bigg|\frac{\partial \theta_i}{\partial u_k}\bigg|\bigg).
\end{equation}
\end{lemma}
\proof
Set $Q=\sinh l_p\sinh l_s\sinh \theta_t$.
By (\ref{cosine law}), we have
\begin{equation*}
\frac{\partial l_i}{\partial r_i}=0,\
\frac{\partial l_i}{\partial r_j}
=\frac{\cosh r_k+\cosh r_j\cosh l_i}{\sinh r_j\sinh l_i},\
\frac{\partial l_i}{\partial r_k}
=\frac{\cosh r_j+\cosh r_k\cosh l_i}{\sinh r_k\sinh l_i}.
\end{equation*}
Since $u_i=\log\tanh\frac{r_i}{2}$, then
\begin{equation*}
\frac{dr_i}{du_i}=\sinh r_i,\ \frac{dr_i}{du_j}=\frac{dr_i}{du_k}=0.
\end{equation*}
Note that
\begin{equation*}
\begin{aligned}
&\cosh l_i
=\cosh \theta_i\sinh l_j\sinh l_k-\cosh l_j\cosh l_k,\\
&\cosh l_j
=\cosh \theta_j\sinh l_i\sinh l_k-\cosh l_i\cosh l_k,\\
&\cosh l_k
=\cosh \theta_k\sinh l_i\sinh l_j-\cosh l_i\cosh l_j.
\end{aligned}
\end{equation*}
Then
\begin{equation*}
\frac{\partial\theta_i}{\partial l_i}=\frac{\sinh l_i}{Q},\
\frac{\partial\theta_i}{\partial l_j}=-\frac{\sinh l_i\cosh \theta_k}{Q},\
\frac{\partial\theta_i}{\partial l_k}=-\frac{\sinh l_i\cosh \theta_j}{Q}.
\end{equation*}
According to the chain rules, we have
\begin{equation}\label{formual 1}
\begin{aligned}
\frac{\partial \theta_i}{\partial u_i}
=&\frac{\partial \theta_i}{\partial l_i}\frac{\partial l_i}{\partial r_i}\frac{dr_i}{du_i}+\frac{\partial \theta_i}{\partial l_j}\frac{\partial l_j}{\partial r_i}\frac{dr_i}{du_i}+\frac{\partial \theta_i}{\partial l_k}\frac{\partial l_k}{\partial r_i}\frac{dr_i}{du_i}\\
=&-\frac{1}{Q}\bigg[\frac{1}{\sinh^2 l_j}\bigg(\cosh l_k\cosh r_k+\cosh l_j\cosh l_k\cosh r_i\\
&+\cosh l_i\cosh l_j\cosh r_k+\cosh l_i\cosh^2 l_j\cosh r_i\bigg)\\
&+\frac{1}{\sinh^2 l_k}\bigg(\cosh l_j\cosh r_j+\cosh l_j\cosh l_k\cosh r_i\\
&+\cosh l_i\cosh l_k\cosh r_j+\cosh l_i\cosh^2 l_k\cosh r_i\bigg)\bigg],
\end{aligned}
\end{equation}
which implies $\frac{\partial \theta_i}{\partial u_i}<0$.
Similarly,
\begin{equation}\label{formual 2}
\begin{aligned}
\frac{\partial \theta_i}{\partial u_j}
=-\frac{1}{Q}\frac{1}{\sinh^2 l_k}\bigg(&-\sinh^2 l_k\cosh r_k+\cosh l_i\cosh r_j+\cosh l_j\cosh r_i\\
&+\cosh l_i\cosh l_k\cosh r_i+\cosh l_j\cosh l_k\cosh r_j\bigg),
\end{aligned}
\end{equation}

\begin{equation}\label{formual 3}
\begin{aligned}
\frac{\partial \theta_i}{\partial u_k}
=-\frac{1}{Q}\frac{1}{\sinh^2 l_j}\bigg(&-\sinh^2 l_j\cosh r_j+\cosh l_i\cosh r_k+\cosh l_k\cosh r_i\\
&+\cosh l_i\cosh l_j\cosh r_i+\cosh l_j\cosh l_k\cosh r_k\bigg).
\end{aligned}
\end{equation}
The formulas (\ref{formual 1}), (\ref{formual 2}) and (\ref{formual 3}) can also be obtained directly from Lemma 4.3 in Guo-Luo \cite{GL2}.
By (\ref{cosine law}), we have
\begin{equation}\label{formual 4}
\begin{aligned}
&\cosh l_i=2\sinh^2\frac{\Phi_{jk}}{2}\sinh r_j\sinh r_k-\frac{1}{2}(e^{r_j-r_k}+e^{r_k-r_j}),\\
&\cosh l_j=2\sinh^2\frac{\Phi_{ik}}{2}\sinh r_i\sinh r_k-\frac{1}{2}(e^{r_i-r_k}+e^{r_k-r_i}),\\
&\cosh l_k=2\sinh^2\frac{\Phi_{ij}}{2}\sinh r_i\sinh r_j-\frac{1}{2}(e^{r_i-r_j}+e^{r_j-r_i}).
\end{aligned}
\end{equation}

We just need to prove that for any positive constants $a,b,c$, the formula (\ref{main formula}) holds, if one of
the following three conditions is satisfied
\begin{description}
\item[(1)] $\lim r_i=+\infty,\ \lim r_j=+\infty,\ \lim r_k=+\infty$,
\item[(2)] $\lim r_i=+\infty,\ \lim r_j=+\infty,\ \lim r_k=c$,
\item[(3)] $\lim r_i=+\infty,\ \lim r_j=a,\ \lim r_k=b$.
\end{description}

For simplicity, we assume $c_i,\ c_j,\ c_k$ and $C_i,\ C_j,\ C_k$ are constants, which are different in the case $\mathbf{(1)}$, $\mathbf{(2)}$ and $\mathbf{(3)}$.

For the case $\mathbf{(1)}$, if $\lim r_i=+\infty,\ \lim r_j=+\infty,\ \lim r_k=+\infty$, then by (\ref{formual 4}), we have
\begin{equation*}
\begin{aligned}
\lim\cosh l_k
=&\lim\left(2\sinh^2\frac{\Phi_{ij}}{2}
\cdot\frac{1}{2}e^{r_i}\cdot\frac{1}{2}e^{r_j}
-\frac{1}{2}(e^{r_i-r_j}+e^{r_j-r_i})\right)\\
=&\lim\left(\frac{1}{2}e^{r_i+r_j}(\sinh^2\frac{\Phi_{ij}}{2}
-e^{-2r_i}-e^{-2r_j})\right)\\
:=&\lim c_ke^{r_i+r_j}.
\end{aligned}
\end{equation*}
Similarly, $\lim\cosh l_i:=\lim c_ie^{r_j+r_k}$ and $\lim\cosh l_j:=\lim c_je^{r_i+r_k}$.
Hence, by (\ref{formual 1}), we have
\begin{equation*}
\begin{aligned}
\lim(-Q\frac{\partial \theta_i}{\partial u_i})
=&\lim\bigg(\frac{1}{c^2_je^{2r_i+2r_k}}\big(c_ke^{r_i+r_j}\cdot\frac{1}{2}e^{r_k}
+c_je^{r_i+r_k}c_ke^{r_i+r_j}\cdot\frac{1}{2}e^{r_i}
+c_ie^{r_j+r_k}c_je^{r_i+r_k}\cdot\frac{1}{2}e^{r_k}\\
&+c_ie^{r_j+r_k}c^2_je^{2r_i+2r_k}\cdot\frac{1}{2}e^{r_i}\big)
+\frac{1}{c^2_ke^{2r_i+2r_j}}\big(c_je^{r_i+r_k}\cdot\frac{1}{2}e^{r_j}
+c_ke^{r_i+r_j}c_je^{r_i+r_k}\cdot\frac{1}{2}e^{r_i}\\
&+c_ie^{r_j+r_k}c_ke^{r_i+r_j}\cdot\frac{1}{2}e^{r_j}
+c_ie^{r_j+r_k}c^2_ke^{2r_i+2r_j}\cdot\frac{1}{2}e^{r_i}\big)\bigg)\\
:=&\lim C_ie^{r_i+r_j+r_k}.
\end{aligned}
\end{equation*}
By (\ref{formual 2}), we have
\begin{equation*}
\begin{aligned}
\lim(-Q\frac{\partial \theta_i}{\partial u_j})
=&\lim\bigg(\frac{1}{c^2_ke^{2r_i+2r_j}}\big(-c^2_ke^{2r_i+2r_j}
\cdot\frac{1}{2}e^{r_k}
+c_ie^{r_j+r_k}\cdot\frac{1}{2}e^{r_j}
+c_je^{r_i+r_k}\cdot\frac{1}{2}e^{r_i}\\
&+c_ie^{r_j+r_k}c_ke^{r_i+r_j}\cdot\frac{1}{2}e^{r_i} +c_je^{r_i+r_k}c_ke^{r_i+r_j}\cdot\frac{1}{2}e^{r_j}\big)\bigg)\\
:=&C_je^{r_k},
\end{aligned}
\end{equation*}
Similarly, by (\ref{formual 3}), we have
\begin{equation*}
\begin{aligned}
\lim(-Q\frac{\partial \theta_i}{\partial u_k})
=&\lim\bigg(\frac{1}{c^2_je^{2r_i+2r_k}}\big(-c^2_je^{2r_i+2r_k}
\cdot\frac{1}{2}e^{r_j}
+c_ie^{r_j+r_k}\cdot\frac{1}{2}e^{r_k}
+c_ke^{r_i+r_j}\cdot\frac{1}{2}e^{r_i}\\
&+c_ie^{r_j+r_k}c_je^{r_i+r_k}\cdot\frac{1}{2}e^{r_i} +c_ke^{r_i+r_j}c_je^{r_i+r_k}\cdot\frac{1}{2}e^{r_k}\big)\bigg)\\
:=&C_ke^{r_j}.
\end{aligned}
\end{equation*}
Since
\begin{align*}
Q&=\sinh l_p\sinh l_s\sinh \theta_t\\
&=\sinh l_p\sinh l_s\sqrt{\cosh^2 \theta_t-1}\\
&=\sqrt{\sinh^2 l_p\sinh^2 l_s\cosh^2 \theta_t-\sinh^2 l_p\sinh^2 l_s}\\
&=\sqrt{(\cosh l_t+\cosh l_p\cosh l_s)^2-\sinh^2 l_p\sinh^2 l_s}\\
&=\sqrt{\cosh^2 l_t+2\cosh l_p\cosh l_s\cosh l_t+\cosh(l_p+l_s)\cosh(l_p-l_s)}\\
&>2,
\end{align*}
then $\lim Q=\lim \sinh l_p\sinh l_s\sinh \theta_t>0$.
For any $C\in (-\infty,+\infty)$, there exists a constant $M>0$, depending on $C$, such that if $r_i\geq M$, then the formula (\ref{main formula}) holds.

For the case $\mathbf{(2)}$, if $\lim r_i=+\infty,\ \lim r_j=+\infty,\ \lim r_k=c$, then by the same calculations, we have $\lim\cosh l_i=\lim c_ie^{r_j}$, $\lim\cosh l_j=\lim c_je^{r_i}$ and $\lim\cosh l_k=\lim c_ke^{r_i+r_j}$.
Furthermore,
\begin{equation*}
\begin{aligned}
\lim(-Q\frac{\partial \theta_i}{\partial u_i})
=\lim C_ie^{r_i+r_j},\
\lim(-Q\frac{\partial \theta_i}{\partial u_j})
=C_j,\
\lim(-Q\frac{\partial \theta_i}{\partial u_k})
=\lim C_ke^{r_j}.
\end{aligned}
\end{equation*}
Note that $\lim Q>0$, for any $C\in (-\infty,+\infty)$, there exists a constant $M>0$, depending on $C$, such that if $r_i\geq M$, then the formula (\ref{main formula}) holds.

For the case $\mathbf{(3)}$: if $\lim r_i=+\infty,\ \lim r_j=a,\ \lim r_k=b$, then by the same calculation,
we have $\lim\cosh l_i=c_i$, $\lim\cosh l_j=\lim c_je^{r_i}$ and $\lim\cosh l_k=\lim c_ke^{r_i}$.
Furthermore,
\begin{equation*}
\begin{aligned}
\lim(-Q\frac{\partial \theta_i}{\partial u_i})
=\lim C_ie^{r_i},\
\lim(-Q\frac{\partial \theta_i}{\partial u_j})
=C_j,\
\lim(-Q\frac{\partial \theta_i}{\partial u_k})
=C_k.
\end{aligned}
\end{equation*}
Note that $\lim Q>0$, for any $C\in (-\infty,+\infty)$, there exists a constant $M>0$, depending on $C$, such that if $r_i\geq M$, then the formula (\ref{main formula}) holds.
\qed

As a direct corollary of Lemma \ref{main lemma}, we have the following result.

\begin{corollary}\label{main corollary}
Suppose $(\Sigma,\mathcal{T})$ is an ideally triangulated surface with boundary.
Let $\Phi: E\rightarrow (0,+\infty)$ be the weight on $(\Sigma,\mathcal{T})$.
Let $n$ be the degree at the boundary component $i\in B$.
Then for any $C_1,C_2,...,C_n\in (-\infty,+\infty)$, there exists a constant $M=M(C_1,C_2,...,C_n)>0$ such that if $r_i\geq M$, then
\begin{equation*}
-\sum_{\{ijk\}\in F}\frac{\partial \theta^{jk}_i}{\partial u_i}
>\sum^n_{j=1,j\sim i}C_j\frac{\partial \theta^{jk}_i}{\partial u_j}.
\end{equation*}
\end{corollary}
\proof
By Lemma \ref{main lemma}, suppose $\{ijk\}\in F$ is a right-angled hyperbolic hexagon adjacent to the boundary components $i,j,k\in B$.
Then for any $C=|C_j|\in (0,+\infty)$, there exists a constant $M=M(C_j)>0$ such that if $r_i\geq M$, then
\begin{equation*}
-\frac{\partial \theta^{jk}_i}{\partial u_i}
=\bigg|\frac{\partial \theta^{jk}_i}{\partial u_i}\bigg|
>C\bigg(\bigg|\frac{\partial \theta^{jk}_i}{\partial u_j}\bigg|+\bigg|\frac{\partial \theta^{jk}_i}{\partial u_k}\bigg|\bigg)
\geq C_j\frac{\partial \theta^{jk}_i}{\partial u_j}.
\end{equation*}
Therefore, for any $C_1,C_2,...,C_n\in (-\infty,+\infty)$, there exists a constant $M=M(C_1,C_2,...,C_n)>0$ such that if $r_i\geq M$, then
\begin{equation*}
-\sum_{\{ijk\}\in F}\frac{\partial \theta^{jk}_i}{\partial u_i}
>\sum_{\{ijk\}\in F}C_j\frac{\partial \theta^{jk}_i}{\partial u_j}
=\sum^n_{j=1,j\sim i}C_j\frac{\partial \theta^{jk}_i}{\partial u_j}.
\end{equation*}
\qed

\begin{lemma}\label{CCF boundary 2}
The solution $u(t)$ of the combinatorial Calabi flow (\ref{Calabi flow formula}) can not reach the boundary $\partial_{0}\Omega(\Phi)$.
\end{lemma}
\proof
Suppose there exists at least $i\in B$ such that $u_i\rightarrow 0^-$,
then $r_i\rightarrow+\infty$ by $u_i=\log\tanh\frac{r_i}{2}$.
By Lemma \ref{Guo-Luo Lemma}, we have $\theta^{jk}_{i}\rightarrow 0$ uniformly as $r_i\rightarrow+\infty$, which implies $K_i\rightarrow 0$ uniformly as $r_i\rightarrow+\infty$.
Thus there exists a number $M_1>0$, such that if $r_i\geq M_1$, $K_i-\overline{K}_i<0$ and $|K_i-\overline{K}_i|\geq\frac{1}{2}\overline{K}_i
\geq\frac{1}{2}M_2>0$.
Note that
\begin{equation}\label{formula 6}
\begin{aligned}
\frac{d u_i}{dt}
&=-\Delta(K-\overline{K})_i\\
&=-\frac{\partial K_i}{\partial u_i}(K_i-\overline{K}_i)-\sum_{j\neq i}\frac{\partial K_i}{\partial u_j}(K_j-\overline{K}_j)\\
&=-\sum_{\{ijk\}\in F}\frac{\partial \theta_i^{jk}}{\partial u_i}(K_i-\overline{K}_i)-\sum_{j\sim i}\left(\frac{\partial \theta^{jk}_i}{\partial u_j}+\frac{\partial \theta^{jl}_i}{\partial u_j}\right)(K_j-\overline{K}_j)\\
&=\bigg[-\sum_{\{ijk\}\in F}\frac{\partial \theta_i^{jk}}{\partial u_i}-\sum_{j\sim i}\left(\frac{\partial \theta^{jk}_i}{\partial u_j}+\frac{\partial \theta^{jl}_i}{\partial u_j}\right)\frac{K_j-\overline{K}_j}{K_i-\overline{K}_i}\bigg](K_i-\overline{K}_i).
\end{aligned}
\end{equation}
By Lemma \ref{CCF lemma}, the combinatorial Calabi energy $\mathcal{C}(u)$ is decreasing along the combinatorial Calabi flow (\ref{Calabi flow formula}),
which implies that for any $i\in B$,
$|K_i-\overline{K}_i|$ is bounded along the combinatorial Calabi flow (\ref{Calabi flow formula}), i.e., there exists a constant $M_3$, such that for any $i\in B$, $|K_i-\overline{K}_i|\leq M_3$.
Therefore, $\big|\frac{K_j-\overline{K}_j}{K_i-\overline{K}_i}\big|
\leq\frac{2M_3}{M_2}$,\ $j\neq i$.
By Corollary \ref{main corollary}, one can choose numbers $C_j=C_j(M_1,M_2,M_3)$, $j=1,2,...,n$, specially $C_j=\frac{4M_3}{M_2}$, if $\frac{\partial\theta^{jk}_i}{\partial u_j}\geq0$, and $C_j=-\frac{4M_3}{M_2}$, if $\frac{\partial\theta^{jk}_i}{\partial u_j}<0$,
such that there exists a number $M_4=M_4(C_1,C_2,...,C_n)>0$, if $r_i\geq \max\{M_1,M_4\}$, then
\begin{equation*}
-\sum_{\{ijk\}\in F}\frac{\partial \theta_i^{jk}}{\partial u_i}-\sum_{j\sim i}\left(\frac{\partial \theta^{jk}_i}{\partial u_j}+\frac{\partial \theta^{jl}_i}{\partial u_j}\right)\frac{K_j-\overline{K}_j}{K_i-\overline{K}_i}
>-\sum_{\{ijk\}\in F}\frac{\partial \theta_i^{jk}}{\partial u_i}-\sum^n_{j=1,j\sim i}C_j\frac{\partial \theta^{jk}_i}{\partial u_j}
>0.
\end{equation*}
By (\ref{formula 6}), we have $\frac{d u_i}{dt}<0.$
The rest of the proof is paralleling to Lemma \ref{Ricci boundary 2}, we omit it for simplicity.
Therefore, the solution of the combinatorial Calabi flow (\ref{Calabi flow formula}) can not reach the boundary $\partial_{0}\Omega(\Phi)$.
\qed

As a direct corollary of Lemma \ref{u(t) bounded 2}, Lemma \ref{CCF boundary 1} and Lemma \ref{CCF boundary 2}, we have the following result on the solution of the combinatorial Calabi flow (\ref{Calabi flow formula}), which implies the longtime existence of the solution $u(t)$ of the combinatorial Calabi flow (\ref{Calabi flow formula}).

\begin{corollary}\label{compact subset 2}
Suppose $(\Sigma,\mathcal{T})$ is an ideally triangulated surface with boundary.
Let $\Phi: E\rightarrow (0,+\infty)$ be the weight on $(\Sigma,\mathcal{T})$.
For any $\overline{K}\in (0,+\infty)^N$ defined on $B=\{1,2,...,N\}$,
the solution $u(t)$ of the combinatorial Calabi flow (\ref{Calabi flow formula}) stays in a compact subset of the admissible space $\Omega(\Phi)$.
As a result, the solution $u(t)$ of the combinatorial Calabi flow (\ref{Calabi flow formula}) exists for all time.
\end{corollary}

The following theorem gives the global convergence of the solution of the combinatorial Calabi flow (\ref{Calabi flow formula}), which is the second part of Theorem \ref{main theorem 1}.

\begin{theorem}
Suppose $(\Sigma,\mathcal{T})$ is an ideally triangulated surface with boundary.
Let $\Phi: E\rightarrow (0,+\infty)$ be the weight on $(\Sigma,\mathcal{T})$.
For any $\overline{K}\in (0,+\infty)^N$ defined on $B=\{1,2,...,N\}$,
the solution of the combinatorial Calabi flow (\ref{Calabi flow formula}) converges exponentially fast.
\end{theorem}
\proof
Note that the solution $u(t)$ of the combinatorial Calabi flow (\ref{Calabi flow formula}) stays in a compact subset of the admissible space $\Omega(\Phi)$ by Corollary \ref{compact subset 2} and the discrete Laplace operator $\Delta=(\frac{\partial K_i}{\partial u_j})_{N\times N}$ is strictly negative by Lemma \ref{matrix negative}.
By the continuity of the eigenvalues of $\Delta$, there exists $\lambda_0>0$ such that the eigenvalues $\lambda_\Delta$ of $\Delta$ satisfies $\lambda_\Delta<-\sqrt{\lambda_0}$ along the combinatorial Calabi flow (\ref{Calabi flow formula}).
Therefore, along the combinatorial Calabi flow (\ref{Calabi flow formula}), we have
\begin{equation*}
\frac{d\mathcal{C}(u(t))}{dt}
=-(K-\overline{K})^T\Delta^2(K-\overline{K})\leq -\lambda_0\mathcal{C}(u(t)),
\end{equation*}
which implies
\begin{equation*}
\mathcal{C}(u(t))
=\frac{1}{2}||K(t)-\overline{K}||^2
\leq e^{-\lambda_0t}||K(0)-\overline{K}||^2.
\end{equation*}
Combining Theorem \ref{rigidity theorem} and Corollary \ref{compact subset 2},
we have
\begin{equation*}
||u(t)-\overline{u}||^2
\leq C_1||K(t)-\overline{K}||^2
\leq C_1e^{-\lambda_0t}||K(0)-\overline{K}||^2
\leq C_2e^{-\lambda_0t}
\end{equation*}
for some positive constants $C_1,C_2$, which completes the proof.
\qed

\begin{remark}
One can also use the Lyapunov Stability Theorem (\cite{Pontryagin}, Chapter 5) to prove the
exponential convergence of the solution $u(t)$ of the combinatorial Calabi flow (\ref{Calabi flow formula}), which is similar to the proof of Theorem \ref{local converge}.
\end{remark}

\section{Fractional combinatorial Calabi flow on surfaces with boundary}\label{section 6}
As the fractional discrete Laplace operator $\Delta^s$ is strictly negative definite on $\Omega(\Phi)$ for any $s\in (-\infty,+\infty)$
by Lemma \ref{matrix negative},
the fractional combinatorial Calabi flow (\ref{FCCF formula}) has many properties similar to that of the combinatorial Calabi flow (\ref{Calabi flow formula}).

\begin{lemma}\label{FCCF lemma}
The function $\mathcal{E}(u)$ defined by (\ref{energy function 1}) and the combinatorial Calabi energy $\mathcal{C}(u)$ defined by (\ref{energy function 2}) are decreasing along the fractional combinatorial Calabi flow (\ref{FCCF formula}).
\end{lemma}

\begin{theorem}\label{local converge for FCCF}
Suppose $(\Sigma,\mathcal{T})$ is an ideally triangulated surface with boundary.
Let $\Phi: E\rightarrow (0,+\infty)$ be the weight on $(\Sigma,\mathcal{T})$.
And $\overline{K}\in (0,+\infty)^N$ is a given function defined on $B=\{1,2,...,N\}$.
If the solution $u(t)$ of the fractional combinatorial Calabi flow (\ref{FCCF formula}) converges to $\overline{u}\in \Omega(\Phi)$, then $K(\overline{u})=\overline{K}$.
Furthermore, for any $\overline{K}\in (0,+\infty)^N$ defined on $B=\{1,2,...,N\}$, there exists a constant $\delta >0$ such that if $||K(u(0))-\overline{K}||=
\sqrt{\sum_{i=1}^{N}(K_i(u(0))-\overline{K}_i)^2}<\delta$, then the solution of the fractional combinatorial Calabi flow (\ref{FCCF formula}) exists for all time and converges exponentially fast.
\end{theorem}

The proof of Lemma \ref{FCCF lemma} is paralleling to that of Lemma \ref{CCF lemma} and
the proof of Theorem \ref{local converge for FCCF} is paralleling to that of Theorem \ref{local converge}.
We omit the details of the proofs here.

Similar to Lemma \ref{u bounded 1} and Lemma \ref{Ricci boundary 1}, we have the following results
on the fractional combinatorial Calabi flow (\ref{FCCF formula}).
As the proofs are almost the same, we also omit the proofs here.

\begin{lemma}\label{u(t) bounded 3}
Suppose $(\Sigma,\mathcal{T})$ is an ideally triangulated surface with boundary.
Let $\Phi: E\rightarrow (0,+\infty)$ be the weight on $(\Sigma,\mathcal{T})$.
For any $\overline{K}\in (0,+\infty)^N$ defined on $B=\{1,2,...,N\}$,
the solution $u(t)$ of the fractional combinatorial Calabi flow (\ref{FCCF formula}) stays in a bounded subset of $(-\infty,0)^N$.
\end{lemma}

\begin{lemma}\label{FCCF boundary 3}
The solution $u(t)$ of the fractional combinatorial Calabi flow (\ref{FCCF formula}) can not reach the boundary $\partial_{l}\Omega(\Phi)$.
\end{lemma}

Lemma \ref{u(t) bounded 3} and Lemma \ref{FCCF boundary 3} show that the solution of the fractional combinatorial Calabi flow (\ref{FCCF formula}) can not reach the boundaries $\partial_{\infty}\Omega(\Phi)$ and $\partial_{l}\Omega(\Phi)$.
However, we do not know how to prove that the solution of the fractional combinatorial Calabi flow (\ref{FCCF formula}) can not reach the boundary $\partial_{0}\Omega(\Phi)$.
As a result, we can not get the longtime existence and global convergence for the solution of the fractional combinatorial Calabi flow (\ref{FCCF formula}) on ideally triangulated surfaces with boundary.
Motivated by the results in \cite{Luo-Xu, Wu-Xu} on fractional combinatorial Calabi flow, we believe the longtime existence and global convergence for the solution of the fractional combinatorial Calabi flow (\ref{FCCF formula}) is still true.
We have the following conjecture.

\begin{conjecture}
Suppose $(\Sigma,\mathcal{T})$ is an ideally triangulated surface with boundary.
Let $\Phi: E\rightarrow (0,+\infty)$ be the weight on $(\Sigma,\mathcal{T})$.
For any $\overline{K}\in (0,+\infty)^N$ defined on $B=\{1,2,...,N\}$, the solution of the fractional combinatorial Calabi flow (\ref{FCCF formula}) exists for all time and converges exponentially fast.
\end{conjecture}
Motivated by \cite{Wu-Xu},
it is believed that a variational formula for the generalized angle  in a right-angled hyperbolic hexagon similar to Glickenstein-Thomas' variational formula for the inner angle in a triangle will play a key role in the proof.

\bigskip

(Xu Xu) School of Mathematics and Statistics, Wuhan University, Wuhan 430072, P.R. China

E-mail: xuxu2@whu.edu.cn\\[2pt]

(Chao Zheng) School of Mathematics and Statistics, Wuhan University, Wuhan 430072, P.R. China

E-mail: czheng@whu.edu.cn\\[2pt]


\begin{thebibliography}{50}
\setlength{\itemsep}{2pt} \small

\bibitem{BPS} A. Bobenko, U. Pinkall, B. Springborn, \emph{Discrete conformal maps and ideal hyperbolic polyhedra}. Geom. Topol. 19 (2015), no. 4, 2155-2215.

\bibitem{Bowers} P. L. Bowers, K. Stephenson, \emph{Uniformizing dessins and Bely\u{i} maps via circle packing}. Mem. Amer. Math. Soc. 170 (2004), no. 805, xii+97 pp.

\bibitem{Chow-Luo} B. Chow, F. Luo, \emph{Combinatorial Ricci flows on surfaces}, J. Differential Geom. 63 (2003), no. 1, 97-129.

\bibitem{Ge1} H. Ge, \emph{Combinatorial methods and geometric equations}, Thesis (Ph.D.)-Peking University, Beijing. 2012. (In Chinese).

\bibitem{Ge2} H. Ge, \emph{Combinatorial Calabi flows on surfaces}, Trans. Amer. Math. Soc. 370 (2018), no. 2, 1377-1391.

\bibitem{Ge3} H. Ge, B. Hua, \emph{On combinatorial Calabi flow with hyperbolic circle patterns}, Adv. Math. 333 (2018), 523-538.

\bibitem{GX3} H. Ge, X. Xu, \emph{$2$-dimensional combinatorial Calabi flow in hyperbolic background geometry}, Differential Geom. Appl. 47 (2016), 86-98.

\bibitem{Ge-Xu 17} H. Ge, X. Xu,  \emph{A discrete Ricci flow on surfaces with hyperbolic background geometry}, Int. Math. Res. Not. IMRN 2017, no. 11, 3510-3527.

\bibitem{Glickenstein} D. Glickenstein,  \emph{Discrete conformal variations and scalar curvature on piecewise flat two and three dimensional manifolds}. J. Differential Geom. 87 (2011), no. 2, 201-237.

\bibitem{GT} D. Glickenstein, J. Thomas,  \emph{Duality structures and discrete conformal variations of piecewise constant curvature surfaces}, Adv. Math. 320 (2017), 250-278.

\bibitem{Gu2} X. D. Gu, R. Guo, F. Luo, J. Sun, T. Wu,  \emph{A discrete uniformization theorem for polyhedral surfaces II},  J. Differential Geom. 109 (2018), no. 3, 431-466.


\bibitem{Gu1} X. D. Gu, F. Luo, J. Sun, T. Wu,  \emph{A discrete uniformization theorem for polyhedral surfaces},  J. Differential Geom. 109 (2018), no. 2, 223-256.

\bibitem{Gu-Luo-Wu} X. D. Gu, F. Luo, T. Wu,  \emph{ Convergence of discrete conformal geometry and computation of uniformization maps}. Asian J. Math. 23 (2019), no. 1, 21-34.

\bibitem{Guo} R. Guo,  \emph{Combinatorial Yamabe flow on hyperbolic surfaces with boundary}. Commun. Contemp. Math. 13 (2011), no. 5, 827-842.


\bibitem{Guo1} R. Guo,  \emph{Local rigidity of inversive distance circle packing}. Trans. Amer. Math. Soc. 363 (2011), no. 9, 4757-4776.

\bibitem{GL2} R. Guo, F. Luo,  \emph{Rigidity of polyhedral surfaces, II}, Geom. Topol. 13 (2009), no. 3, 1265-1312.

\bibitem{Li-Xu-Zhou} S.Y. Li, X. Xu, Z. Zhou, \emph{Combinatorial Yamabe flow on hyperbolic bordered surfaces}. \href{https://arxiv.org/abs/2204.08191}{ arXiv:2204.08191v1[math.GT]}.

\bibitem{Luo1} F. Luo, \emph{Combinatorial Yamabe flow on surfaces}, Commun. Contemp. Math. 6 (2004), no. 5, 765-780.

\bibitem{Luo3} F. Luo, \emph{Rigidity of polyhedral surfaces, III}. Geom. Topol. 15 (2011), no. 4, 2299-2319.

\bibitem{Luo-Wu} F. Luo, T. Wu, \emph{Koebe conjecture and the Weyl problem for convex surfaces in hyperbolic $3$-space}. \href{https://arxiv.org/abs/1910.08001}{ arXiv:1910.08001v2[math.GT]}.

\bibitem{Luo-Xu} Y. Luo, X. Xu, \emph{Combinatorial Calabi flows on surfaces with boundary}. Calc. Var. Partial Differential Equations 61 (2022), no. 3, Paper No. 81, 12 pp.

\bibitem{Pontryagin} L.S. Pontryagin, \emph{Ordinary differential equations}, Addison-Wesley Publishing Company Inc., Reading, 1962.


\bibitem{Ratcliffe} J.G. Ratcliffe, \emph{Foundations of hyperbolic manifolds}. Second edition. Graduate Texts in Mathematics, 149, xii+779 pp. Springer, New York (2006). ISBN: 978-0387-33197-3; 0-387-33197-2.


\bibitem{Stephenson} K. Stephenson, \emph{Introduction to Circle Packing}. The Theory of Discrete Analytic Functions. Cambridge University Press, Cambridge (2005).

\bibitem{Sun-Wu-Gu-Luo} J. Sun, T. Wu, X. D. Gu, F. Luo, \emph{Discrete conformal deformation: algorithm and experiments}. SIAM J. Imaging Sci. 8 (2015), no. 3, 1421-1456.

\bibitem{Thurston} W. Thurston, \emph{Geometry and topology of $3$-manifolds}, Princeton lecture notes, 1976.

\bibitem{Wu} T. Wu, \emph{Finiteness of switches in discrete Ricci flow}. Master Thesis, Tsinghua University, Beijing, 2014. (In Chinese).

\bibitem{Wu-Gu-Sun} T. Wu, X. D. Gu, J. Sun, \emph{Rigidity of infinite hexagonal triangulation of the plane}. Trans. Amer. Math. Soc. 367 (2015), no. 9, 6539-6555.

\bibitem{Wu-Xu} T. Wu, X. Xu, \emph{Fractional combinatorial Calabi flow on surfaces}, \href{https://arxiv.org/abs/2107.14102}{ arXiv:2107.14102[math.GT]}.

\bibitem{Wu-Zhu} T. Wu, X. Zhu, \emph{The convergence of discrete uniformizations for closed surfaces}. \href{https://arxiv.org/abs/2008.06744}{ arXiv:2008.06744v2[math.GT]}. To appear in J. Differential Geom.

\bibitem{Xu18} X. Xu, \emph{Rigidity of inversive distance circle packings revisited}. Adv. Math. 332 (2018), 476-509.

\bibitem{Xu MRL} X. Xu, \emph{A new proof of Bowers-Stephenson conjecture}, Math. Res. Lett. 28 (2021), no. 4, 1283-1306.

\bibitem{Xu22} X. Xu, \emph{A new class of discrete conformal structures on surfaces with boundary}. Calc. Var. Partial Differential Equations 61 (2022), no. 4, Paper No. 141, 23 pp.
    
\bibitem{Xu CAG} X. Xu, \emph{Parameterized discrete uniformization theorems and curvature flows for polyhedral surfaces, I},  \href{https://arxiv.org/abs/1806.04516} {arXiv:1806.04516 [math.GT].} To appear in Comm. Anal. Geom.

\bibitem{Xu1} X. Xu, \emph{Rigidity and deformation of discrete conformal structures on polyhedral surfaces}.
    \href{https://arxiv.org/abs/2103.05272}{ arXiv:2103.05272v2[math.GT]}.
    
\bibitem{XZ TAMS} X. Xu, C. Zheng, \emph{Parameterized discrete uniformization theorems and curvature flows for polyhedral surfaces, II},
Trans. Amer. Math. Soc. 375 (2022), no. 4, 2763-2788.

\bibitem{XZ CVPDE} X. Xu, C. Zheng, \emph{Prescribing discrete Gaussian curvature on polyhedral surfaces}. Calc. Var. Partial Differential Equations 61 (2022), no. 3, Paper No. 80, 17 pp.

\bibitem{Zhang-Guo-Zeng-Luo-Yau-Gu} M. Zhang, R. Guo, W. Zeng, F. Luo, S.T. Yau, X. D. Gu,  \emph{The unified discrete surface Ricci flow}. Graphical Models., 76(5), 321-339 (2014).

\bibitem{Zhu-Xu} X. Zhu, X. Xu, \emph{Combinatorial Calabi flow with surgery on surfaces}, Calc. Var. Partial Differential Equations 58 (2019), no. 6, Paper No. 195, 20 pp.




\end{thebibliography}
\end{document}